\documentclass[11pt]{article}  
\usepackage{geometry}                		
\usepackage{graphicx}
\usepackage{amsmath}
\usepackage{epstopdf}
\usepackage{mathptmx}
\usepackage{graphicx}
\usepackage{graphics}
\usepackage{subfig}
\usepackage{caption}
\usepackage{epsfig}
\usepackage{url}
\usepackage{color}

\begin{document}

\title{Reduced temporal convergence rates in high-order splitting schemes}

\author{M. T. Warnez\footnote{Engineering Physics, University of Michigan, College of Engineering, 2355 Bonisteel Blvd., Ann Arbor, MI 48109-2104, USA, mwarnez@umich.edu} \and B. K. Muite\footnote{Computer, Electrical and Mathematical Sciences and Engineering Division, King Abdullah University of Science \& Technology (KAUST), P.O. Box 4700, Thuwal 23955, Saudi Arabia, benson.muite@kaust.edu.sa \newline \textcolor{white}{\,\,\,\,\,\,\,\,\,\,\,}\textit{Keywords}: splitting methods, high-order accuracy, reaction-diffusion equations \newline \textcolor{white}{\,\,\,\,\,\,\,\,\,\,\,}\textit{Mathematics Subject Classification}: 65M12, 35K57, 65M70}}

%

\maketitle

\begin{abstract}
Recently-derived high-order splitting schemes with complex coefficients are shown to exhibit reduced convergence rates for certain parabolic evolution equations.  When applied to semilinear reaction-diffusion equations with periodic boundary conditions, these splitting schemes are most efficient when the diffusion terms are much less stiff than the reaction terms.  An explanation for this is given in a simple setting.
\end{abstract}

\section{Introduction} \label{intro}

Splitting methods constitute an effective and widely-used class of numerical integrators \cite{HLW05,LR04}.  While splitting schemes with real-valued splitting coefficients are limited to second order accuracy \cite{GK96,S91}, recent work \cite{CCDV09,HO09} has shown that higher order methods can be constructed by using complex-valued coefficients.  Motivated by reaction-diffusion problems, optimized coefficients for splitting methods of orders 6 and 8 were developed in \cite{BCCM13}.  Although their theoretical framework applies only to linear parabolic systems, their numerical results indicate that nonlinear scalar equations can be treated as well.  No attempt to solve either linear or nonlinear vector problems was reported.  To clarify the scope of these optimized splitting schemes, this work investigates splitting solutions to the diffusion equation with a potential and the Gray-Scott reaction-diffusion system. The findings reported here seem also to hold for a wide variety of other reaction diffusion equations, including the Fisher, Lotka-Volterra and Zeldovich equations. It should also be noted that there is a large literature on solving reaction diffusion equations, and the Gray-Scott equations in particular \cite{EO13,K04,KT05,P93,ZWZ08}. We do not aim at completeness in reviewing this literature, nor claim that the splitting methods used here are the most efficient for all of these problems. The aim of this work is to show that the optimized coefficients reported in \cite{BCCM13} are most effective for moderate time steps when the diffusion constant is small.  For a large value of the diffusion constant, or for highly oscillatory solutions, these \emph{optimized} coefficients are not optimal. In such cases it may be better to derive equation-dependent, and possibly even solution-dependent, splitting coefficients.

\section{Summary of Previous Work}\label{SummaryPrevious}

Reaction-diffusion equations for a reactive species $u(\mathbf{x},t)$ take the general form
\begin{align}
\frac{\partial u}{\partial t} = D\Delta u + R(u)
\end{align}
where $D$ is the diffusivity constant and $R(u)$ is some nonlinear function of $u$. In splitting methods, two sets of equations are solved in an alternating manner: a linear differential equation, such as the heat equation
\begin{align}
\frac{\partial u}{\partial t} = D\Delta u 
\end{align}
and a set of ordinary differential equations
\begin{align}
\frac{\partial u}{\partial t} = R(u).
\end{align}
To obtain high-order methods, \cite{BCCM13} used composition techniques with complex coefficients. The underlying integrators chosen in  \cite{BCCM13} are Lie-Trotter and Strang splittings.  Given the evolution equation
\begin{align}
\frac{\partial u}{\partial t} = A(u) + B(u)
\end{align}
with operators $A$ and $B$, the Lie-Trotter approximation is
\begin{align}
u(\delta t)\approx \exp(At)\exp(Bt)u(0).
\end{align}
A higher order approximation is obtained by
\begin{align}
u(\delta t)\approx \prod_{i=1}^N\exp(c(i)At)\exp(d(i)Bt)u(0)
\end{align}
where $c(i)$ and $d(i)$ are complex coefficients chosen to minimize some functional related to the error in the approximation and order conditions of the scheme. For the Strang splitting method
\begin{align}
u(\delta t)\approx \exp(At/2)\exp(Bt)\exp(At/2)u(0)
\end{align}
a higher order approximation is obtained by
\begin{align}
u(\delta t)\approx  \prod_{i=1}^N\exp(\gamma(i)At/2)\exp(\gamma(i)Bt)\exp(\gamma(i)At/2)u(0)
\end{align}
where $\gamma(i)$ are complex coefficients chosen to minimize some functional related to the error.

\section{Linear Diffusion Equation with a Potential}

We consider a modification of the linear diffusion equation considered in \cite{BCCM13}, which takes the form
\begin{align} \label{PL}
\frac{\partial u}{\partial t} = Du_{xx} + (3+\sin(10x))\cos(12x)u.
\end{align}
For this we adopt the splitting strategy of \cite{BCCM13}, splitting the heat equation
\begin{align} \label{PLeq}
\frac{\partial u}{\partial t} = Du_{xx}
\end{align}
from the ordinary differential equations
\begin{align}
\frac{\partial u}{\partial t} = (3+\sin(10x))\cos(12x)u.
\end{align}
In the Fourier space, denoted by $\,\,\widehat{\,}\,\,$, the linear equation has the solution
\begin{align} \label{PLsol}
\widehat{u} = \widehat{u_0}e^{-Dk^2t}
\end{align}
where $k$ is the wavenumber. The ordinary differential equations are solved by
\begin{align}
u = u_0\exp\left[t(3+\sin(10x))\cos(12x)\right].
\end{align}

\setcounter{subfigure}{0}
\begin{figure}
	\subfloat[Convergence with high diffusivity]{\includegraphics[width=.3\paperwidth]{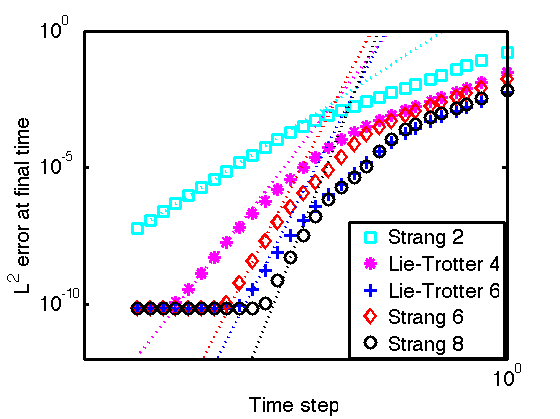}}
	\subfloat[Convergence with low diffusivity]{\includegraphics[width=.3\paperwidth]{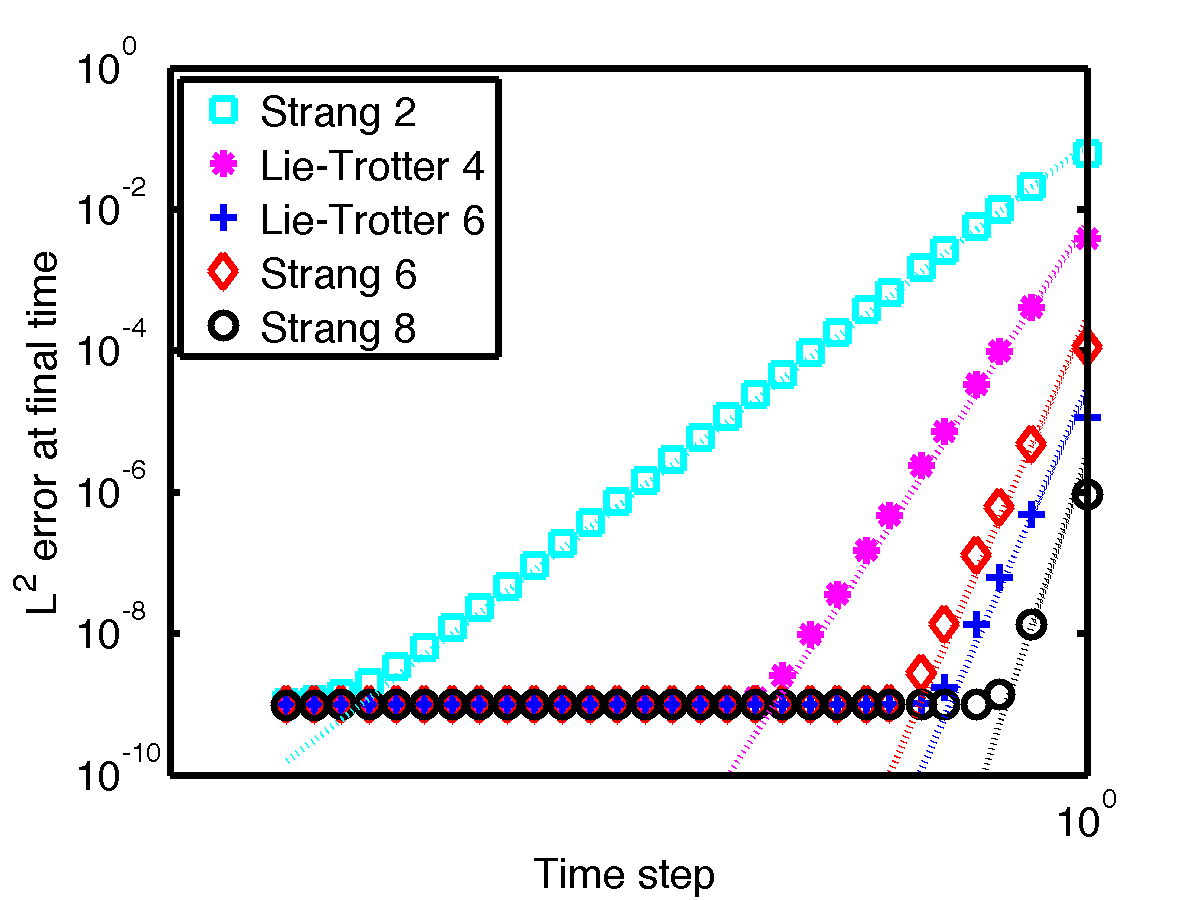}}\\
	\subfloat[Computational efficiency with high diffusivity]{\includegraphics[width=.3\paperwidth]{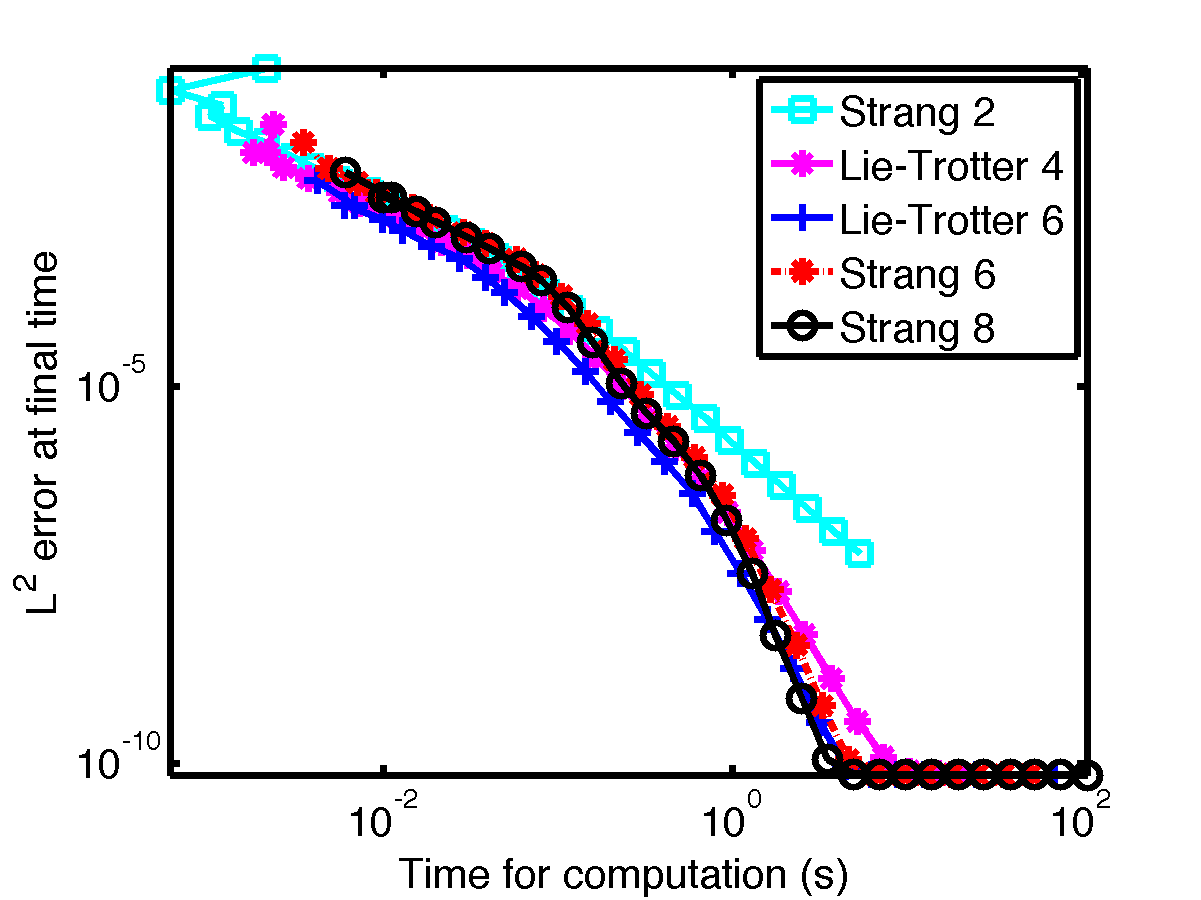}}
	\subfloat[Computational efficiency with low diffusivity]{\includegraphics[width=.3\paperwidth]{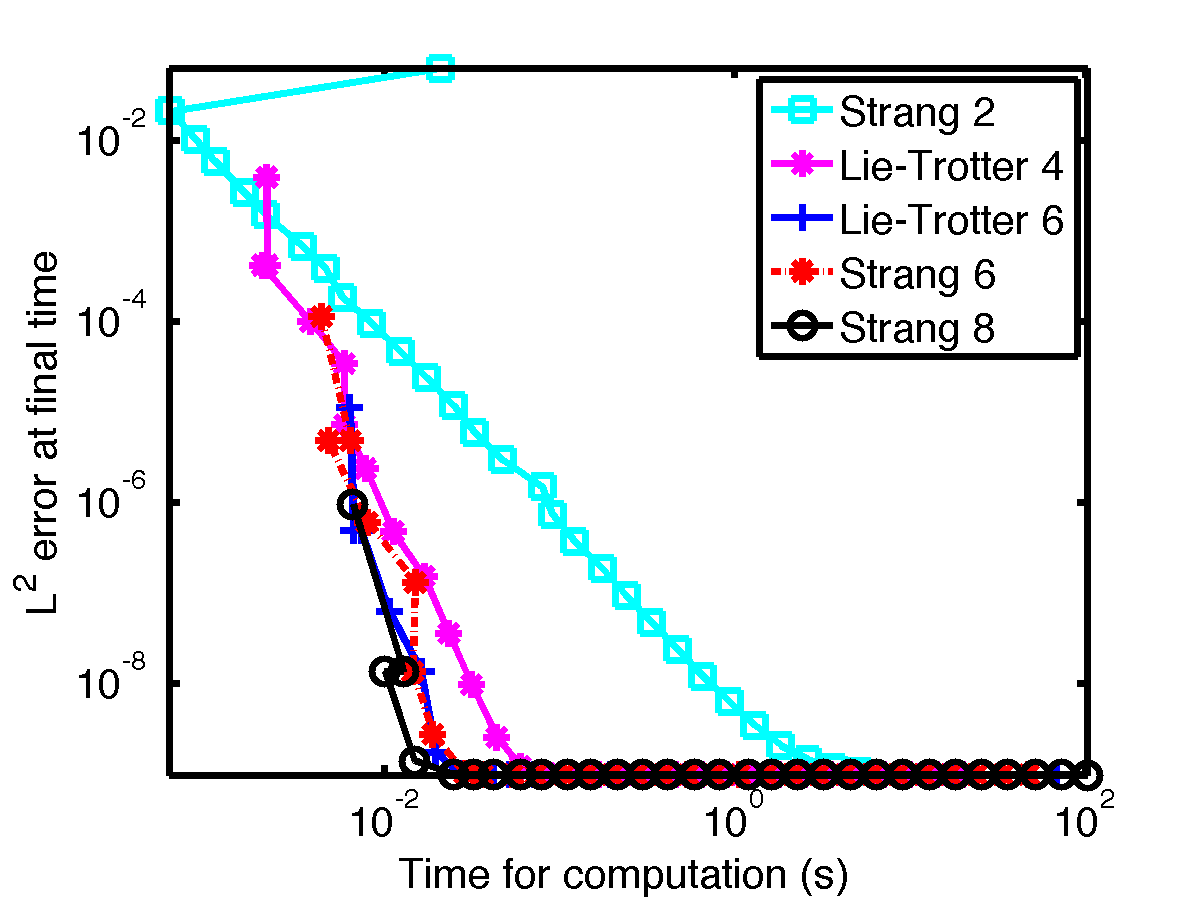}}
	\caption{Numerical results for a linear diffusion equation with a potential, eq.\ \eqref{PL}. The absolute discrete $l^2$ error is shown, where $l^2:=\sqrt{\sum_i \lvert u_{\mathrm{ref},i}-u_i\rvert^2}\Delta x$, and the reference solution $u_\mathrm{ref}$ is computed with a time step of $2^{-20}$ using the eighth order composition method. The solid lines in (a) and (b) show expected convergence rates for second, fourth, sixth and eighth order methods. Results were obtained from a splitting scheme in which the heat equation is used as the second operator. In the high-diffusion case $D=10$ and in the low-diffusion case $D=0.005$. The initial condition was $u=0.5-0.2\exp(\sin(8x))$. A total of 1024 grid points uniformly spaced over $x\in[-\pi,\pi)$ were used.}
	\label{fig:LinearPotential2}
\end{figure}	

Results of numerical experiments performed in MATLAB R2013a on a 2.53GHz Intel i5 dual core chip are in Fig.\ \ref{fig:LinearPotential2}. The implementation can be found in the supplementary materials. The expected order of convergence is only obtained when the product of the time step and the diffusion constant is relatively small. The figure also shows that the highest order methods are not always the most efficient. For these particular examples, the sixth order Lie-Trotter splitting method is computationally most efficient for moderate accuracies.

\section{Gray-Scott Equations}

The Gray-Scott equations describe an autocatalytic reaction-diffusion process between two chemical constituents with concentrations $u$ and $v$:
\begin{align}
\frac{\partial u}{\partial t} &= D_u u_{xx} + \alpha\left(1-u\right) - uv^2,\\
\frac{\partial v}{\partial t} &= D_v v_{xx} - \beta v + uv^2.
\end{align}
The chemicals have independent diffusion constants $D_u$ and $D_v$.  The parameter $\alpha > 0$ determines the rate at which the reactant is added, while $\beta >0$ determines the rate at which the product is subtracted.  By splitting the linear terms of the reactive contribution into the linear system,
\begin{align} \label{GSuLeq}
\frac{\partial u}{\partial t} &= D_u u_{xx} + \alpha\left(1-u\right), \\
\frac{\partial v}{\partial t} &= D_v v_{xx} - \beta v, \label{GSvLeq}
\end{align}
the corresponding nonlinear system
\begin{align}
\frac{\partial u}{\partial t} &= -uv^2 \\
\frac{\partial v}{\partial t} &=  uv^2
\end{align}
has the conservation relation $u + v = u_0 + v_0$.  This allows the nonlinear system to be uncoupled to obtain $v_t=(u_0+v_0-v)v^2$, which can be integrated, see for example \cite{G12}, to obtain,
\begin{align}
u &= u_0 + v_0 - v \\
v &= \frac{u_0 + v_0}{1 + W_0\left[\frac{u_0}{v_0}\exp\left(\frac{u_0}{v_0}-\left(u_0 + v_0\right)^2t\right)\right]} \label{GSvNsol}
\end{align}
where $W_0$ is the principal branch of the Lambert $W$ function \cite{CGHJK96,C13,V12}.  Writing eq.\ \eqref{GSvNsol} in terms of $W_0$ is only possible when $v_0 \ne 0$. In most physical applications $v_0,u_0\geq0$, although alternative representations may be more appropriate when $v_0$ is close to zero.  It should also be noted that the Lambert $W$ function can be used to solve a wide variety of ordinary differential equations that appear in modeling chemical reactions \cite{CGHJK96,C13,V12}, and is useful for splitting schemes applied to other reaction diffusion equations, such as for the Zeldovich equation, $u_t=D_uu_{xx}+u(1-u)(u-\alpha)$.

In the Fourier space, the linear system eq.\ \eqref{GSuLeq}-eq.\ \eqref{GSvLeq} has the solution
\begin{align}
\widehat{u} &= \left(\widehat{u_0} - \frac{\widehat{\alpha}}{\alpha + D_uk^2}\right)e^{-\left(\alpha + D_uk^2\right)t} + \frac{\widehat{\alpha}}{\alpha + D_uk^2}, \\
\widehat{v} &= \widehat{v_0}e^{-\left(\beta + D_vk^2\right)t}.
\end{align}

\setcounter{subfigure}{0}
\begin{figure}
	\subfloat[Convergence with high diffusivity]{\includegraphics[width=.3\paperwidth]{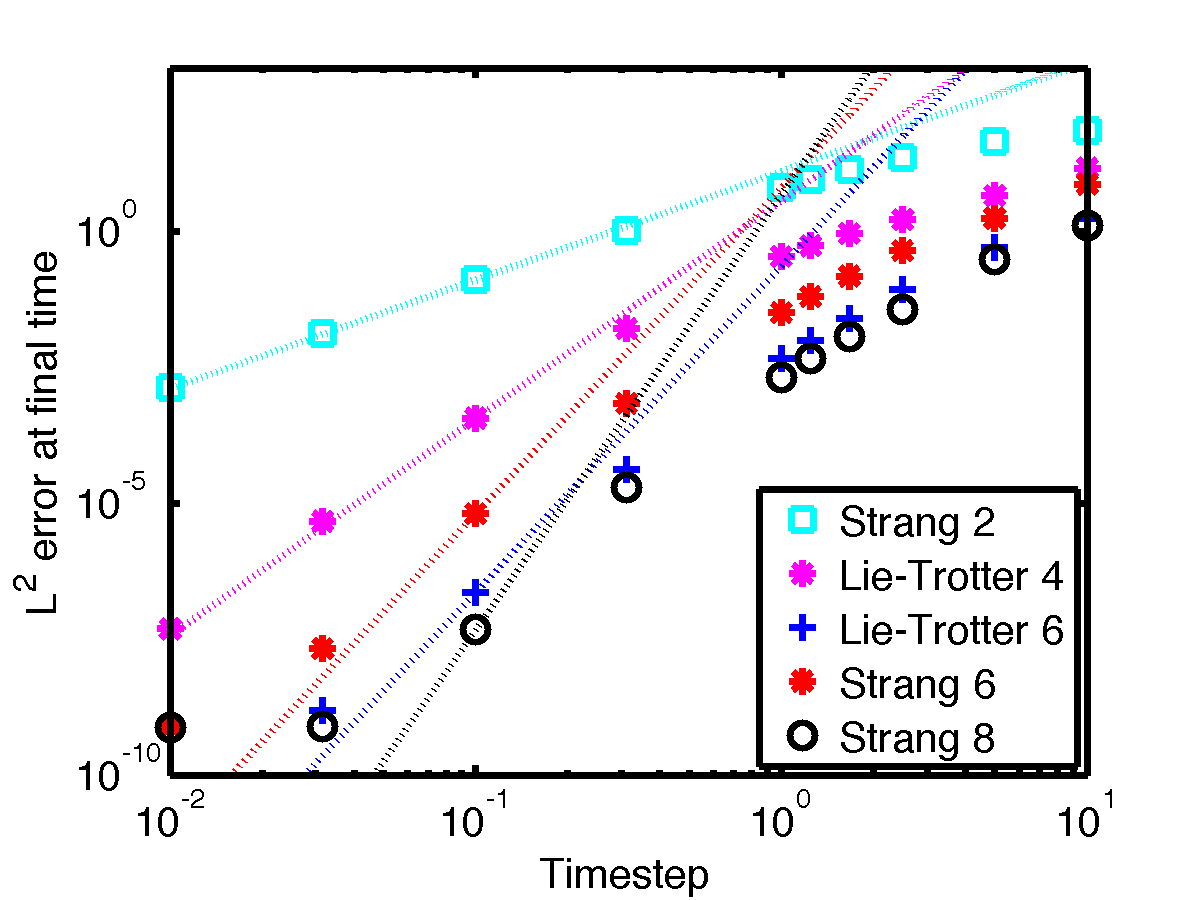}}
	\subfloat[Convergence with low diffusivity]{\includegraphics[width=.3\paperwidth]{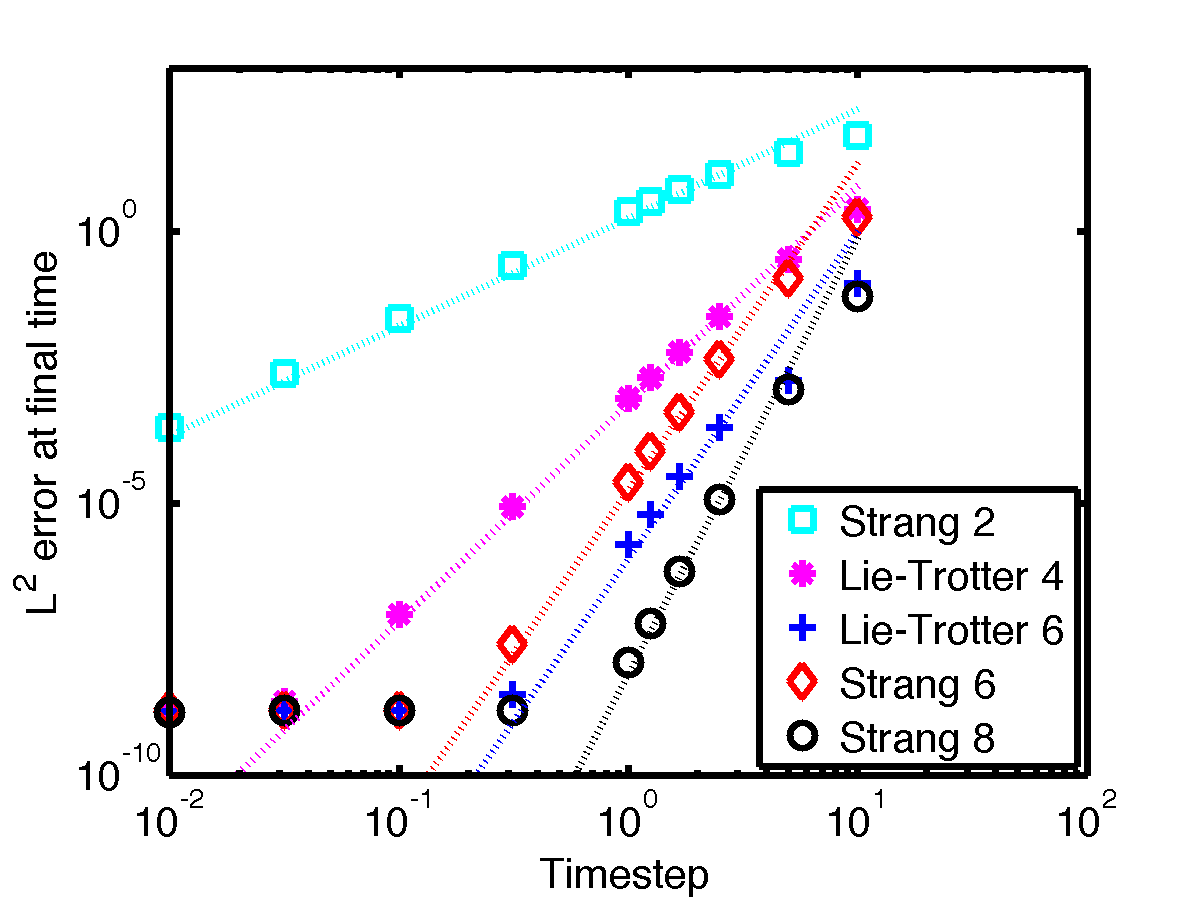}}\\
	\subfloat[Computational efficiency with high diffusivity]{\includegraphics[width=.3\paperwidth]{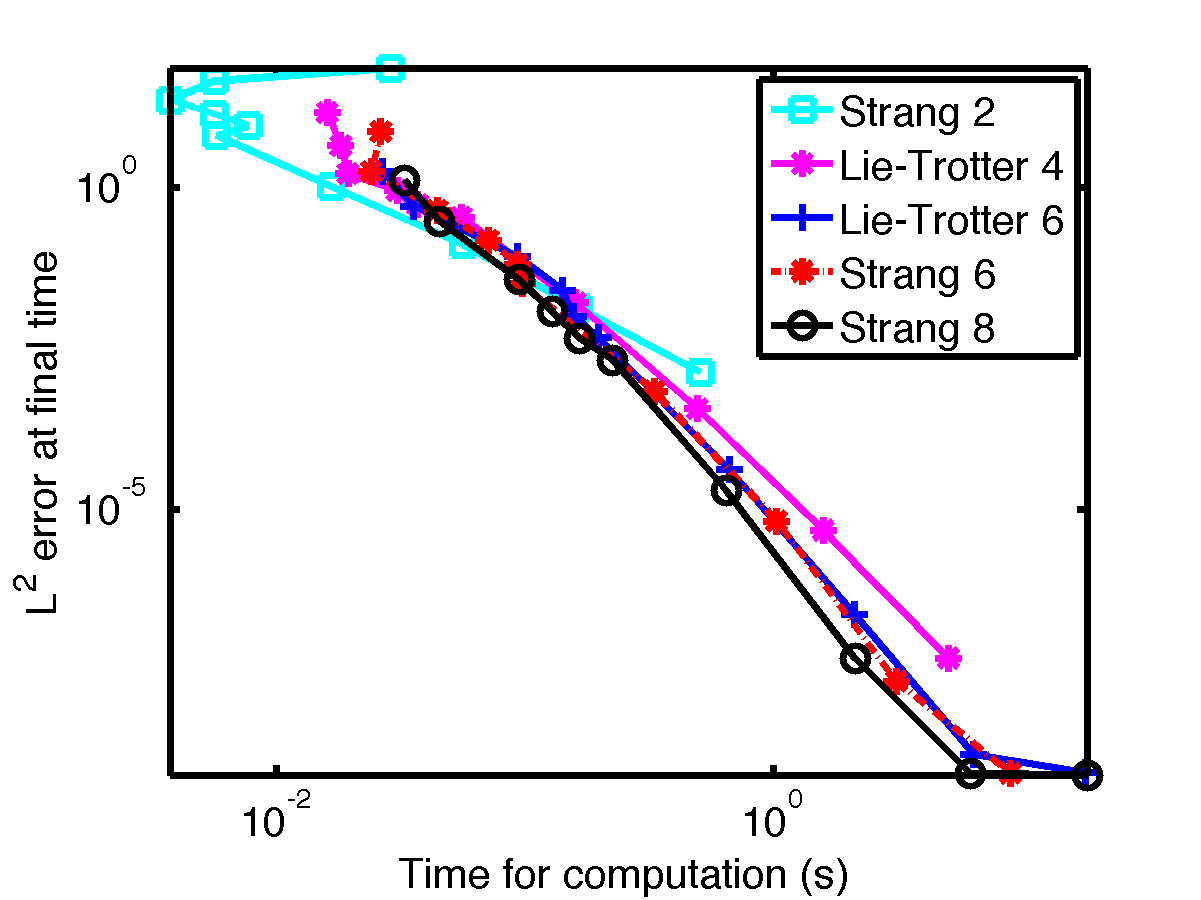}}
	\subfloat[Computational efficiency with low diffusivity]{\includegraphics[width=.3\paperwidth]{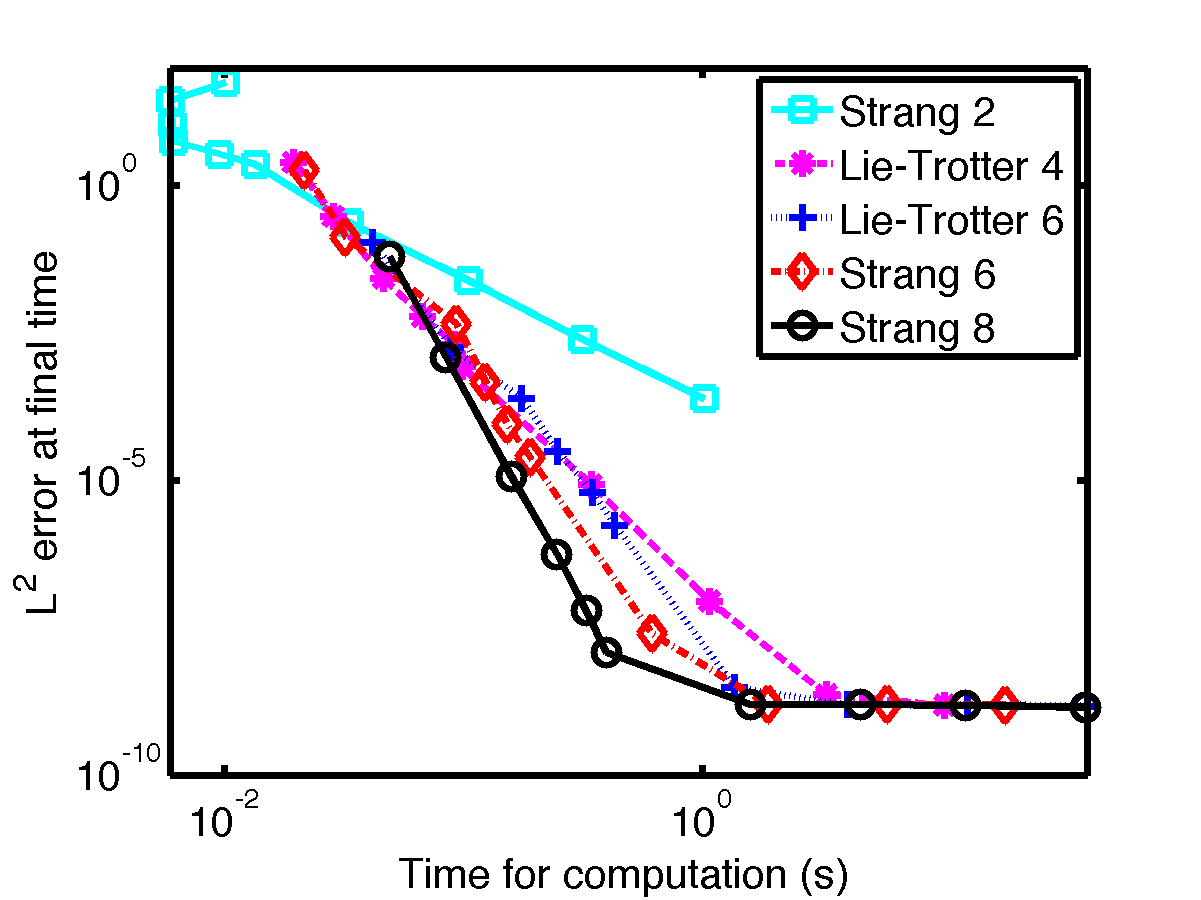}}
	\caption{Numerical results for the Gray-Scott equations. The absolute error, defined as $\lVert u_\mathrm{ref}-u\rVert_2+\lVert v_\mathrm{ref}-v\rVert_2$, is shown. The reference solutions were computed with a time step of $10^{-4}$ using the eighth order composition method. In the low diffusivity case $D_u=D_v=0.001$, while in the high diffusivity case $D_u=1$ and $D_v=0.01$. The dotted lines in (a) and (b) show expected convergence rates for second, fourth, sixth and eighth order methods. Results were obtained with a splitting scheme where the linear equation is used as the first operator. For both cases, the parameters were $\alpha=0.09$ and $\beta=0.086$ and the simulation was run until a final time of 10. The initial conditions were $u=0.0903(1+0.9\exp(\cos(3x)))$, $v=0.952(1-0.9\sin(\cos(x)))$. A total of 512 uniform grid points uniformly spaced over $x\in[-\pi,\pi)$ were used.}
	\label{fig:GrayScott}
\end{figure}

Figure \ref{fig:GrayScott} contains results obtained in MATLAB R2013a on a 2.53GHz Intel i5 dual core chip, which show that the expected order of convergence is obtained for low, but not high, diffusion constants. The figure also shows that the differences in efficiency for the fourth, sixth and eight order methods are small in the high diffusion case, but higher order methods become very efficient in the low diffusion case. For these particular examples, the eighth order Strang splitting method is computationally most efficient. In place of MATLAB's symbolic Lambert $W$ function, a fully numerical implementation from \cite{G13} was used to reduce computation times.  The implementations can be found in the supplementary materials. 

\subsection{Example Application}

The Gray-Scott system is known to exhibit many interesting behaviors (such as chaotic solutions, non-constant stationary solutions, and a variety of dynamic pattern formation regimes) which are best explored by numerical simulation or through a combination of analysis and numerical simulation \cite{DKZ97,KBHS09,MR04,MPSS96,NU01,PW09,S06,WW03}.  Numerous online reports \cite{AACHNS13,L08,M13} include extensive two dimensional simulations of the Gray Scott equations, but a full numerical investigation of the behavior of solutions to the Gray-Scott equations in three dimensions has not been performed. To demonstrate that the high accuracy splitting schemes can be used when accurate solutions are required, we use a splitting scheme to demonstrate the self-reproducing dynamics described in \cite{NU01} and the chaotic dynamics described in \cite{KBHS09}. Results of the simulations are shown in figures \ref{fig:1Dbreakup} and \ref{fig:1Dchaos}. The implementations can be found in the supplementary materials. In these simulations, $v$ can be close to zero, for which it becomes problematic to numerically evaluate the exact solution in eq.\ \eqref{GSvNsol}. Consequently, within the splitting scheme, the implicit midpoint rule is used to solve $v_t=(u_0+v_0-v)v^2$ instead of the exact nonlinear solution. As explained in \cite{BCCM13}, this technique also gives a high-order splitting method. The uniformly valid  implicit midpoint rule has been chosen here for demonstration purposes, though it may be possible to find good analytical approximations of the exact solution that are fast to evaluate for small values of $v$.

\begin{figure}
\center
\includegraphics[width=.6\paperwidth]{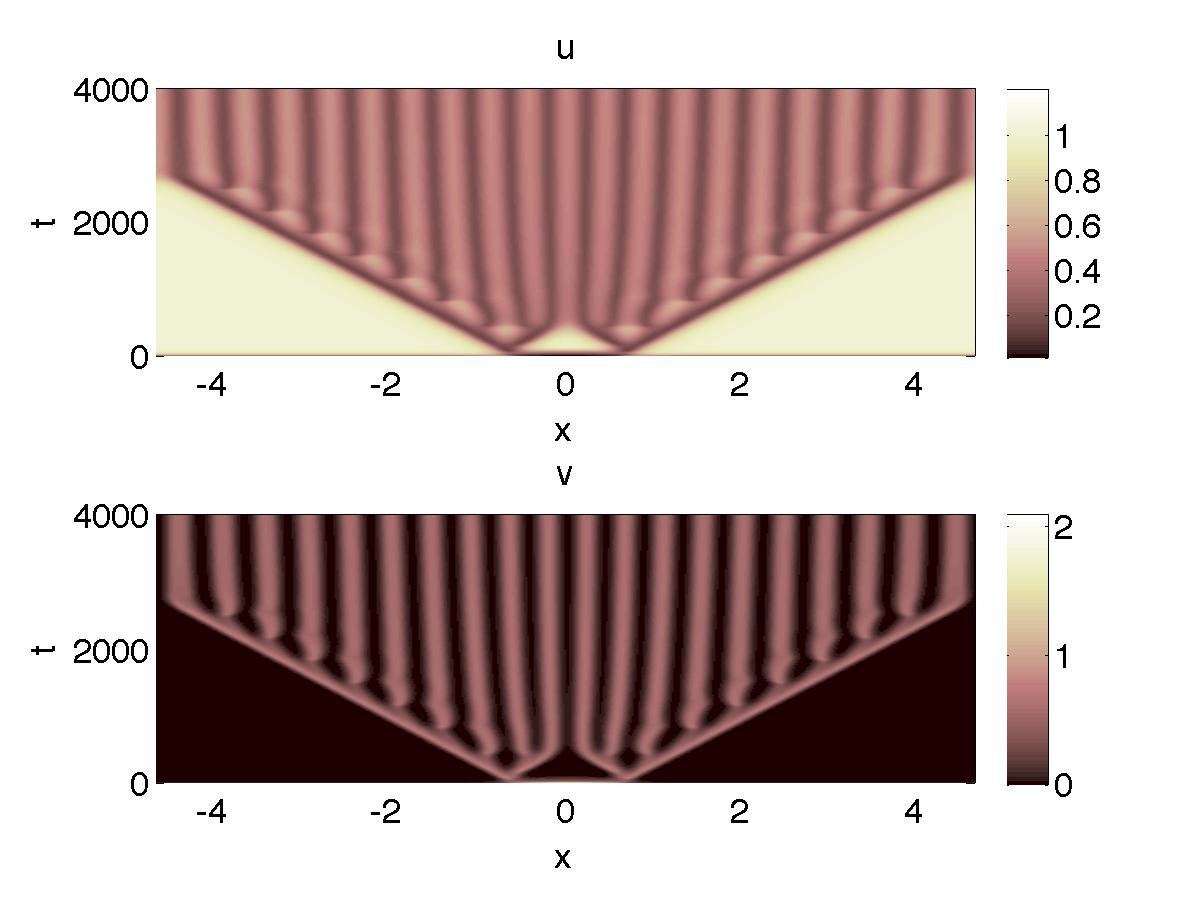}
\caption{The results of a numerical simulation of the Gray-Scott equations on $[-1.5\pi,1.5\pi)$ with $\alpha=0.04$, $\beta=0.1$, $D_u=0.001$, $D_v=0.0001$, $u(t=0,x)= \exp(-2x^2)$, $v(t=0,x) = 0.1+ \exp(-4x^2)$, a time step of 1 and 512 grid points.}\label{fig:1Dbreakup}
\end{figure}

\begin{figure}
\center
\includegraphics[width=.6\paperwidth]{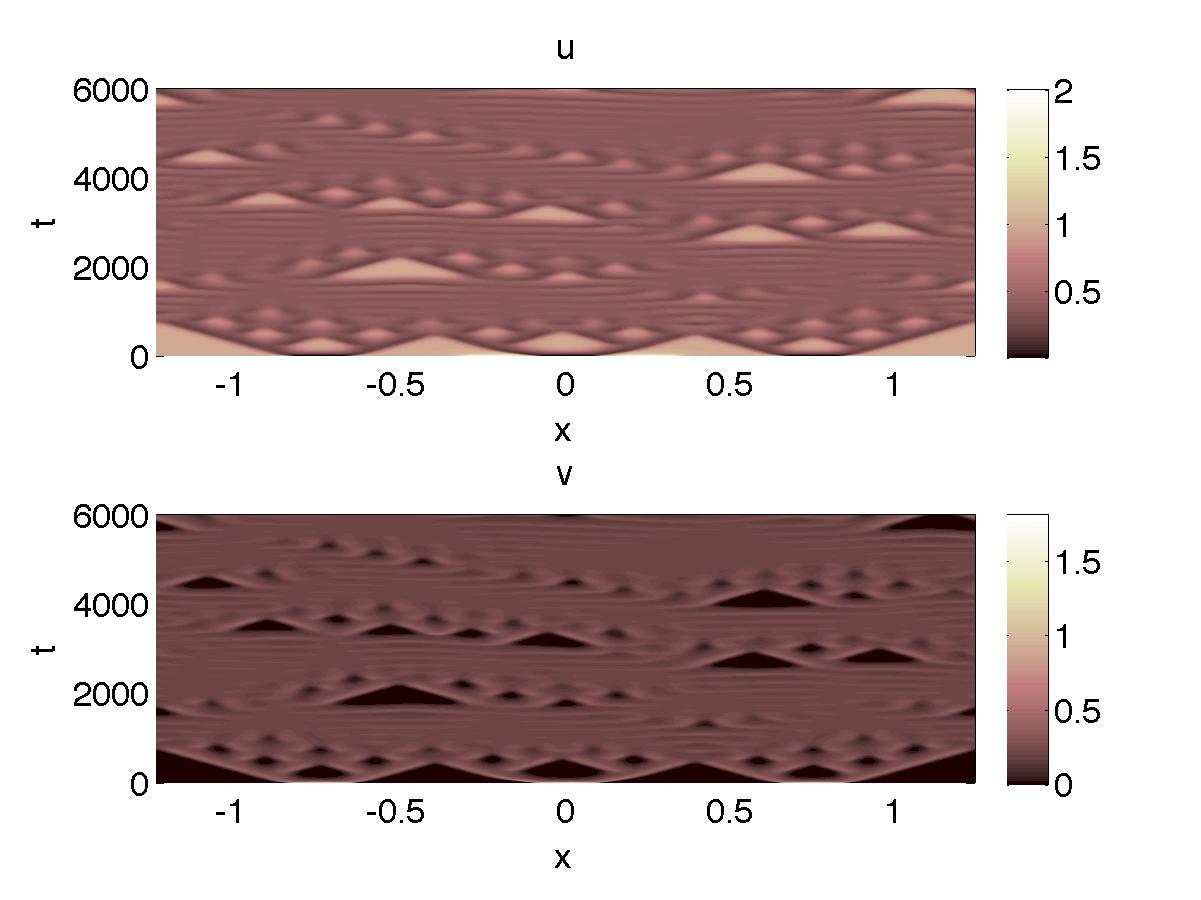}
\caption{The results of a numerical simulation of the Gray-Scott equations on $[-1.25,1.25)$ with $\alpha=0.028$, $\beta=0.081$, $D_u=2\times10^{-5}$, $D_v=10^{-5}$, $u(t=0,x)=1+\exp(-2(10(x-0.25))^8)+\exp(-2(10(x+0.23))^8)$, $v(t=0,x) =\exp(-4(10x)^6)+\exp(-4(10(x-0.75))^6)$, a time step of 0.25 and 256 grid points.}\label{fig:1Dchaos}
\end{figure}

\section{Explanation of Reduced Convergence Rates}

To understand why large diffusion constants lead to reduced convergence rates, we examine the nature of the splitting error for the linear diffusion equation with a potential. The analysis is similar to that in \cite{S00}. Consider 
\begin{align}
u_t=D_uu_{xx}+(3+\sin(10x))\cos(12x)u
\end{align}
with $A:=D_u\partial_{xx}$ and $B:=(3+\sin(10x))\cos(12x)$. The solution approximated by one time step of Strang splitting is
\begin{align}
u(\delta t)=\exp\left(B\delta t/2\right)\exp\left(A\delta t\right)\exp\left(B\delta t/2\right)u(0)
\end{align}
while the exact solution is given by
\begin{align}
u(\delta t)=\exp\left[(A+B)\delta t\right]u(0).
\end{align}
A Taylor expansion of the exponential function gives the error after one time step as
\begin{align}
\hspace{2em}&{}\frac{\delta t^3}{3!}\left[- BAB - \frac{BAA}{4} + \frac{BBA}{2} + \frac{ABB}{2} - \frac{AAB}{4} - \frac{ABA}{2} \right]u+O((\delta t)^4) \notag
\\&{}=\frac{\delta t^3}{3!}\left[D\frac{((3+\sin(10x))\cos(12x))^2u_{xx}+\left[((3+\sin(10x))\cos(12x))^2u\right]_{xx} }{4} \right. \notag
\\&{}\hspace{3em} - D^2\frac{(3+\sin(10x))\cos(12x)u_{xxxx}+\left[(3+\sin(10x))\cos(12x)u\right]_{xxxx}}{2} \notag
\\&{}\hspace{3em}+D^2[(3+\sin(10x))\cos(12x)u_{xx}]_{xx} \notag
\\&{}\hspace{3em}+D(3+\sin(10x))\cos(12x)[(3+\sin(10x))\cos(12x)u]_{xx} \bigg] +O((\delta t)^4). \label{eq:errorexpansion}
\end{align}
Table \ref{Tab:ErrorMagnitude} shows that if $D$ is not small, the coefficients of the the $O((\delta t)^3)$ error term can be large, leading to lower than expected accuracy if the time step is not particularly small. For higher order schemes, the first error term in the Taylor approximation contains a larger number of terms, with greater powers of the derivatives, for which reduced convergence rates are easier to observe if the time step is not sufficiently small. A similar effect occurs in the nonlinear case.
 
\begin{table}[htdp]
\begin{center}
\begin{tabular}{|c|c|c|}
\hline
Term & $D=10$ & $D=0.01$ \\
\hline
$D(3+\sin(10x))\cos(12x)[(3+\sin(10x))\cos(12x)u]_{xx}$  &  $1.2\times10^4$ & $1.2\times10^{1}$\\
$D^2(3+\sin(10x))\cos(12x)u_{xxxx}/2$            & $6.7\times10^5$ & $6.7\times10^{-1}$\\
$D((3+\sin(10x))\cos(12x))^2u_{xx}/2$             & $1.8\times10^{3}$ & $1.8\,\,\,\,\,$\\
$D\left[((3+\sin(10x))\cos(12x))^2u\right]_{xx}/2$   & $5.7\times10^3$ & $5.7\,\,\,\,\,$\\
$D^2\left[(3+\sin(10x))\cos(12x)u\right]_{xxxx}/2$  & $1.1\times10^6$ & $1.1\,\,\,\,\,$\\
$D^2[(3+\sin(10x))\cos(12x)u_{xx}]_{xx}$    & $2.7\times10^5$ &$2.7\times10^2\,\,\,\,\,$\\
\hline
\end{tabular}
\end{center}
\caption{The size of the terms in eq.\ \eqref{eq:errorexpansion} in the $l^{\infty}$ norm for the initial condition $u=0.5-0.2\exp(\sin(8x))$ used in generating Fig.\ \ref{fig:LinearPotential2}.}
\label{Tab:ErrorMagnitude}
\end{table}%
  
\section{Conclusions}

Reduced convergence rates for high-order splitting methods with complex coefficients for parabolic equations seem not to have been encountered before.  Sufficiently small times steps are needed in order for the first component of the error term to be significantly smaller than the computed components of the Taylor expansion that are a good approximation of the exact solution. These reduced convergence rates are most apparent when the diffusion term is large, and can make higher-order splitting methods significantly less efficient for solving partial differential equations.  In some cases, therefore, the use of low order splitting methods may be more appropriate.  

Further examination of the error terms may suggest better ways to generate optimized splitting schemes for semilinear parabolic partial differential equations. In particular, the methods used in \cite{BCCM13} do not minimize error terms for a specific differential equation. The splitting constants are given, but the methods used to find these are not. The results in Table \ref{Tab:ErrorMagnitude}, demonstrate that certain error terms are more important than others, therefore, for some applications, weighting the error terms in a splitting coefficient optimization scheme may give significantly more accurate results.
 
Similar splitting errors are observed for splitting methods applied to nonlinear Schr\"{o}dinger equations in the small dispersion limit \cite{BJM02,DT10,RS13}, for which the nonlinear term produces large errors. The reduced convergence rates observed here resemble observations of order reduction for Runge-Kutta methods applied to parabolic partial differential equations with periodic boundary conditions (see \cite{HO05,KR11}), as well as order reduction due to the effects of boundary conditions (see \cite{FOS12,HV03}). While the mechanism of reduced convergence is different, the effect on the numerical scheme is comparable. Finally, similar splitting errors are also observed in splitting/fractional step methods for boundary value problems \cite[p.~258]{B00}. In these, a grid space-dependent approximation error may have a large coefficient, so that a good approximation can only be obtained for very small grid sizes.

\section{Acknowledgements}
The authors thank Winfried Auzinger, John Boyd, Yiannis Hadjimichael, Lajos Loczi, David Ketcheson, Christian Klein, Paul Rigge, and Kristelle Roidot for helpful discussions. A portion of this work was completed at KAUST, to whom MTW thanks for support through a summer internship.  BKM thanks the University of Michigan for support during the initial phase of this project.

\end{document}